\documentclass[a4paper,10pt]{amsart}

\usepackage{amsmath}
\usepackage{amssymb}
\usepackage{dsfont}
\usepackage{enumerate}
\usepackage{eurosym}
\usepackage[all]{xy}

\newcommand{\0}{\emptyset}

\renewcommand{\a}{\alpha}
\renewcommand{\b}{\beta}
\newcommand{\g}{\gamma}
\renewcommand{\l}{\lambda}

\newcommand{\q}{\mathds{Q}}

\newcommand{\z}{\mathds{Z}}
\newcommand{\n}{\mathds{N}}

\newcommand{\m}{{\mathfrak{M}}}

\newtheorem*{main-theorem-A}{Theorem A}
\newtheorem*{main-theorem-B}{Theorem B}
\newtheorem*{question}{Question}
\newtheorem*{answer}{Answer}
\newtheorem*{proposition*}{Proposition}

\newtheorem{theorem}{Theorem}[section]
\newtheorem{proposition}[theorem]{Proposition}

\theoremstyle{definition}
\newtheorem*{main-conjecture}{Conjecture}

\newcommand{\opname}[1]{\operatorname{\mathsf{#1}}}

\newcommand{\Hom}{\opname{Hom}}
\newcommand{\Ext}{\opname{Ext}}

\newcommand{\End}{\opname{End}}

\newcommand{\rep}{\opname{rep}}

\newcommand{\udim}{\underline{\dim}\,}

\begin{document}

\title[A remark on real root representations of quivers]{A remark on the constructibility of real root representations of quivers using universal extension functors}

\author{Marcel Wiedemann}

\address{Department of Pure Mathematics, University of Leeds, Leeds LS2 9JT, U.K.}

\email{marcel@maths.leeds.ac.uk}

\date{\today}


\begin{abstract}
In this paper we consider the following question: Is it possible to construct all real root representations of a given quiver $Q$ by using universal extension functors, starting with a real Schur representation? We give a concrete example answering this question negatively.
\end{abstract}

\maketitle

\setcounter{section}{-1}
\section{Introduction}
\label{introduction}
\noindent Let $k$ be a field and let $Q$ be a (finite) quiver. We fix a representation $S$ with $\End_{kQ} S=k$ and $\Ext^1_{kQ} (S,S)=0$. In analogy to \cite[Section 1]{ringel} we consider the following subcategories of $\rep_k Q$. Let $\m^S$ be the full subcategory of all modules $X$ with $\Ext^1_{kQ}(S,X)=0$ such that, in addition, $X$ has no direct summand which can be embedded into some direct sum of copies of $S$. Similarly, let $\m_S$ be the full subcategory of all modules $X$ with $\Ext^1_{kQ}(X,S)=0$ such that, in addition, no direct summand of $X$ is a quotient of a direct sum of copies of $S$. Finally, let $\m^{-S}$ be the full subcategory of all modules $X$ with $\Hom_{kQ} (X,S)=0$, and let $\m_{-S}$ be the full subcategory of all modules $X$ with $\Hom_{kQ} (S,X)=0$. Moreover, we consider
\begin{eqnarray*}
    \m_S^S = \m^S \cap \m_S,\quad   \m_{-S}^{-S}= \m^{-S} \cap \m_{-S}.
\end{eqnarray*}
According to \cite[Proposition 1 \& $1^*$ and Proposition 2]{ringel}, we have the following equivalences of categories
\begin{eqnarray*}
\overline{\sigma}_S&:& \m^{-S} \to \m^S/S,		\\
\underline{\sigma}_S&:& \m_{-S} \to \m_S/S,		\\
\sigma_S&:& \m_{-S}^{-S} \to \m^S_S/S,
\end{eqnarray*}
where $\m^{S}/S$ denotes the quotient category of $\m^{S}$ modulo the maps which factor through direct sums of copies of $S$, similarly for $\m_S/S$ and $\m^{S}_S/S$. We call the functor $\sigma_S$ {\it universal extension functor}. A brief description of these functors is given in Section \ref{notation}.
This paper is dedicated to the following question.
\begin{question}[$\star$]
Let $\a$ be a positive non-Schur real root for $Q$ and let $X_\a$ be the unique indecomposable representation of dimension vector $\a$. 

Does there exist a sequence of real Schur roots $\b_1,\ldots,\b_n\; (n\ge 2)$ such that
\begin{eqnarray*}
			X_\a = \sigma_{X_{\b_n}}\cdot \ldots \cdot \sigma_{X_{\b_2}} (X_{\b_1})\quad ?
\end{eqnarray*}
Here, $X_{\b_i}$ denotes the unique indecomposable representation of dimension vector $\b_i$.
\end{question}
One might reformulate the above question as follows. Is it possible to construct all real root representations of $Q$ using universal extension functors, starting with a real Schur representation?

One of the nice facts about the universal extension functor $\sigma_S$ is that it allows one to keep track of certain properties of representations. For instance, the functor $\sigma_S$ preserves indecomposable tree representations \cite[Lemma 3.16]{wiedemann} (for a definition of ``tree representation'' and background results we refer the reader to \cite[Introduction]{ringel_3}) and, moreover, if we apply the functor $\sigma_S$ to a representation of known endomorphism ring dimension, we can easily compute the dimension of the endomorphism ring of the resulting representation \cite[Proposition 3 \& $3^*$]{ringel}. Hence, if $X_\a=\sigma_{X_{\b_n}}\cdot \ldots \cdot \sigma_{X_{\b_2}} (X_{\b_1})$ with $\b_i\; (i=1,\ldots,n)$ real Schur roots, then $X_\a$ is a tree representation and one can easily compute $\dim \End_{kQ} X_\a$.

Question $(\star )$ was first answered affirmatively by Ringel \cite[Section 2]{ringel} for the quiver
\begin{center}
$Q(g,h): 	\xymatrix{  1 \ar@/^1.8pc/[r]^{\mu_1}
				         \ar@/^0.9pc/[r]^*-<0.4pc>{\vdots}_{\mu_g} & 2 \ar@/^1.8pc/[l]^{\nu_h} \ar@/^0.9pc/[l]^*-<1.3pc>{\vdots}_{\nu_1}}$,
\end{center}
with $g,h\ge 1$. In \cite[Theorem B]{wiedemann} Question $(\star)$ was answered affirmatively for the quiver
\begin{center}

$Q(f,g,h)$: 		$\xymatrix{ 1 \ar@<1.5ex>[r]^{\l_1} \ar@<-1ex>[r]^*-<0.4pc>{\vdots}_{\l_f} & 2 \ar@/^1.8pc/[r]^{\mu_1}
				         \ar@/^0.9pc/[r]^*-<0.4pc>{\vdots}_{\mu_g} & 3 \ar@/^1.8pc/[l]^{\nu_h} \ar@/^0.9pc/[l]^*-<1.3pc>{\vdots}_{\nu_1}}$,

\end{center}
with $f,g,h\ge 1$. More examples of real root representations which can be constructed using universal extension functors can be found in \cite[Appendix]{wiedemann-thesis}.

Hence, there are quivers for which Question $(\star)$ can be answered affirmatively. The question is, can it be answered affirmatively in general? Unfortunately the answer is negative in general.

\begin{answer}[to Question $(\star)$]
In Section \ref{counterexample} we give a concrete example answering Question $(\star )$ negatively.
\end{answer}

This paper is organized as follows. In Section \ref{notation} we discuss further notation and background results and in Section \ref{counterexample} we describe an example answering Question $(\star )$ negatively.

{\bf{Acknowledgements.}} The author would like to thank his supervisor, Prof. W. Crawley-Boevey, for his continuing support and guidance. The author also wishes to thank Prof. C. Ringel for his interest in this work and for stimulating discussions.

\section{Further Notation and Background Results}
\label{notation}

Let $k$ be a field. Let $Q$ be a finite quiver, i.e. an oriented graph with finite vertex set $Q_0$ and finite arrow set $Q_1$ together with two functions $h,t:Q_1 \to Q_0$ assigning head and tail to each arrow $a\in Q_1$. A representation $X$ of $Q$ is given by a vector space $X_i$ (over $k$) for each vertex $i\in Q_0$ together with a linear map $X_a: X_{t(a)} \to X_{h(a)}$ for each arrow $a\in Q_1$. Let $X$ and $Y$ be two representations of $Q$. A homomorphism $\phi: X \to Y$ is given by linear maps $\phi_i: X_i\to Y_i$ such that for each arrow $a\in Q_1$, $a:i\to j$ say, the square
\begin{center}
$\xymatrix{ X_i \ar@<0.0ex>[r]^-{X_a} \ar@<0.0ex>[d]_{\phi_i} & X_j \ar@<0.0ex>[d]^-{\phi_j}	\\
					  Y_i \ar@<0.0ex>[r]^-{Y_a} & Y_j }$
\end{center}
commutes. 

A dimension vector for $Q$ is given by an element of $\n^{Q_0}$. We will write $e_i$ for the coordinate vector at vertex $i$ and by $\a[i],\, i\in Q_0,$ we denote the $i$-th coordinate of $\a\in \n^{Q_0}$. We can partially order $\n^{Q_0}$ via $\a\ge \b$ if $\a[i]\ge \b[i]$ for all $i\in Q_0$. We define $\a >\b$ to mean $\a\ge \b$ and $\a\ne \b$. If $X$ is a finite dimensional representation, meaning that all vector spaces $X_i\; (i\in Q_0)$ are finite dimensional, then $\udim X= (\dim X_i)_{i\in Q_0}$ is the dimension vector of $X$. Throughout this paper we only consider finite dimensional representations. We denote by $\rep_k Q$ the full subcategory with objects the finite dimensional representations of $Q$.
The Ringel form on $\z^{Q_0}$ is defined by
\begin{eqnarray*}
    \langle\a,\b\rangle = \sum_{i\in Q_0} \a[i]\b[i]   -   \sum_{a\in Q_1}\a[t(a)]\b[h(a)]
\end{eqnarray*}
Moreover, let $(\a,\b)=\langle\a,\b\rangle+\langle\b,\a\rangle$ be its symmetrization. 

We say that a vertex $i\in Q_0$ is loop-free if there are no arrows $a:i\to i$. By a quiver without loops we mean a quiver with only loop-free vertices. For a loop-free vertex $i\in Q_0$ the simple reflection $s_i:\z^{Q_0}\to \z^{Q_0}$ is defined by
\begin{eqnarray*}
	s_i(\a):=\a-(\a,e_i)e_i.
\end{eqnarray*}

A simple root is a vector $e_i$ for $i\in Q_0$. The set of simple roots is denoted by $\Pi$. The Weyl group, denoted by $W$, is the subgroup of $\textrm{GL}(\z^n)$, where $n=|Q_0|$, generated by the $s_i$. By $\Delta^+_{\textrm{re}}(Q) := \{\a \in W(\Pi) : \a> 0\} $ we denote the set of (positive) real roots for $Q$. 

We have the following remarkable theorem.

\begin{theorem}[{Kac \cite[Theorem 1 and 2]{kac}, Schofield \cite[Theorem 9]{schofield}}]
\label{kac-theorem}
Let $k$ be a field, $Q$ be a quiver and let $\a \in \Delta^+_{\textrm{re}}(Q)$. There exists a unique indecomposable representation (up to isomorphism) of dimension vector $\a$.
\end{theorem}

For finite fields and algebraically closed fields the theorem is due to Kac \cite[Theorem 1 and 2]{kac}. As pointed out in the introduction of \cite{schofield}, Kac's method of proof showed that the above theorem holds for fields of characteristic $p$. The proof for fields of characteristic zero is due to Schofield \cite[Theorem 9]{schofield}.

For a given positve real root $\a$ for $Q$ the unique indecomposable representation (up to isomorphism) of dimension vector $\a$ is denoted by $X_\a$. By a real root representation we mean an $X_\a$ for $\a$ a positive real root. A Schur representation is a representation with $\End_{kQ} (X)=k$. By a real Schur representation we mean a real representation which is also a Schur representation. A positive real root is called a real Schur root if $X_\a$ is a real Schur representation.

We have the following useful formula: if $X,Y$ are representations of $Q$ then we have 
\begin{eqnarray*}
	\dim \Hom_{kQ}(X,Y) - \dim \Ext^1_{kQ} (X,Y) = \langle \udim X, \udim Y \rangle.
\end{eqnarray*}
It follows that $\Ext^1_{kQ}(X_\a,X_\a)=0$ for $\a$ a real Schur root.

\subsection{Universal Extension Functors}
\label{universal}

We use this section to describe briefly how the functors
\begin{eqnarray*}
\overline{\sigma}_S&:& \m^{-S} \to \m^S/S,		\\
\underline{\sigma}_S&:& \m_{-S} \to \m_S/S,		\\
\sigma_S&:& \m_{-S}^{-S} \to \m^S_S/S,
\end{eqnarray*}
operate on objects.

The functor $\overline{\sigma}_S$ is given by the following construction: Let $X\in\m^{-S}$ and let $E_1,\ldots,E_r$ be a basis of the $k$-vector space $\Ext^1_{kQ} (S,X)$. Consider the exact sequence $E$ given by the elements $E_1,\ldots,E_r$
\begin{eqnarray*}
		E: 0\to X\to Z\to \bigoplus_r S \to 0.
\end{eqnarray*}
According to \cite[Lemma 3]{ringel} we have $Z\in \m^S$ and we define $\overline{\sigma}_S(X):=Z$. Now, let $Y\in\m_{-S}$ and let $E'_1,\ldots,E'_s$ be a basis of the $k$-vector space $\Ext^1_{kQ}(Y,S)$. Consider the exact sequence $E'$ given by $E'_1,\ldots,E'_s$
\begin{eqnarray*}
		E': 0\to \bigoplus_s S\to U\to Y\to 0.
\end{eqnarray*}
Then we have $U\in \m_S$ and we set $\underline{\sigma}_S(Y):=U$. The functor $\sigma_S$ is given by applying both constructions successively.

The inverse $\overline{\sigma}_S^{-1}$ is constructed as follows: Let $X\in \m^{S}$ and let $\phi_1,\ldots,\phi_r$ be a basis of the $k$-vector space $\Hom_{kQ}(X,S)$. Then by \cite[Lemma 2]{ringel} the sequence
\begin{eqnarray*}
		0\to X^{-S}\to X \stackrel{(\phi_i)_i}{\longrightarrow} \bigoplus_r S\to 0
\end{eqnarray*}
is exact, where $X^{-S}$ denotes the intersection of the kernels of all maps $X\to S$. We set $\overline{\sigma}^{-1}_S (X):=X^{-S}$. Now, let $Y\in\m_{S}$. The inverse $\underline{\sigma}_S^{-1}$ is given by $\underline{\sigma}_S^{-1}(Y):=Y/Y'$, where $Y'$ is the sum of the images of all maps $S\to Y$. The inverse $\sigma^{-1}_S$ is given by applying both constructions successively.

Both constructions show that
\begin{equation}
 		\udim \sigma^{\pm 1}_S(X) = \udim X - (\udim X, \udim S)\cdot\udim S. \tag{$\dagger$}\label{refl-dim}
\end{equation}

Moreover, we have the following proposition.

\begin{proposition}[{\cite[Proposition 3 \& $3^*$]{ringel}}]
Let $X\in \m^{-S}_{-S}$. Then
\begin{eqnarray*}
	\dim \End_{kQ} \sigma_S(X) = \dim \End_{kQ}(X) + \langle \udim X, \udim S \rangle \cdot \langle \udim S,\udim X \rangle.
\end{eqnarray*}
Let $Y\in \m^{S}_{S}$. Then
\begin{eqnarray*}
	\dim \End_{kQ} \sigma^{-1}_S(Y) = \dim \End_{kQ}(Y) - \langle \udim Y, \udim S \rangle \cdot \langle \udim S,\udim Y \rangle.
\end{eqnarray*}
\end{proposition}

\section{A negative and unpleasant example}
\label{counterexample}

Let $k$ be a field and let $Q$ be a quiver. We recall Question $(\star )$ stated in the introduction.
\begin{question}[$\star$]
Let $\a$ be a positive non-Schur real root for $Q$ and let $X_\a$ be the unique indecomposable representation of dimension vector $\a$.

Does there exist a sequence of real Schur roots $\b_1,\ldots,\b_n\; (n\ge 2)$ such that
\begin{eqnarray*}
			X_\a = \sigma_{X_{\b_n}}\cdot \ldots \cdot\sigma_{X_{\b_2}} (X_{\b_1})\quad ?
\end{eqnarray*}
\end{question}
\noindent We remark that in the case that $X_\a$ can be constructed in the above way we have $\b_i<\a$ for $i=1,\ldots,n$.

In the following we give an explicit example of a non-Schur real root representations which cannot be constructed using universal extension functors.

\pagebreak

We consider the quiver $Q$

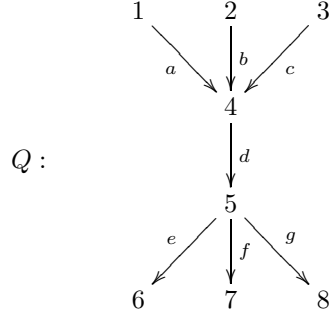
\begin{picture}(2,4.5)(-4.5,-4) 
\put(-1.5,-2){$Q:$}
$   \xymatrix{ 1 \ar@<.0ex>[dr]_a                &2  \ar@<.0ex>[d]^b                      	   &3  \ar@<.0ex>[dl]^c   \\
																								 &4	 \ar@<.0ex>[d]^d									 			   &   \\
																								 &5	 \ar@<.0ex>[dl]_e \ar@<.0ex>[d]^f \ar@<.0ex>[dr]^g  & 	 \\
						   6                                 &7                                            &8 }$
\end{picture}

\noindent and the real root $\a=(1,1,1,8,12,2,7,7)=s_8s_7s_5s_4s_8s_7s_5s_8s_7s_5s_6s_4s_5s_4s_1s_2s_3(e_4)$.

For the convenience of the reader we give an explicit description of the representation $X_\a$. 

We start by considering the representation $X_\a$ over the field $k=\q$. In this case, one can use the result \cite[Proposition A.4]{crawley} to construct the representation $X_\a$; we get

\begin{picture}(2,5)(-2.2,-4) 
\put(0.5,-2){$X_\a:$}
$   \xymatrix{     		&&k  \ar@<.0ex>[dr]_{X_a}             &k  	  \ar@<.0ex>[d]^{X_b}    		                  	     &k  \ar@<.0ex>[dl]^{X_c}   \\
											&&																	  &k^8	  \ar@<.0ex>[d]^{X_d}												 			   &   \\
											&&																	  &k^{12}	\ar@<.0ex>[dl]_{X_e} \ar@<.0ex>[d]^{X_f} \ar@<.0ex>[dr]^{X_g}  & 	 \\
						        	&& k^2                                &k^7                                     				           &k^7 }$
\end{picture}
\bigskip

\noindent with

\begin{eqnarray*}
	X_a	&=& \left[ \begin{array}{cccccccc}
									0	& 0 & 0 & 0 & 0 & 1 & 0 & 0
								 \end{array}
					\right]^t, \\
			& &	\\
  X_b	&=& \left[ \begin{array}{cccccccc}
									0	& 0 & 0 & 0 & 0 & 0 & 1 & 0
								 \end{array}
					\right]^t, \\
			& &	\\		
	X_c	&=& \left[ \begin{array}{cccccccc}
									0	& 0 & 0 & 0 & 0 & 1 & 1 & 1
								 \end{array}
					\right]^t,\\
			& & \\
	X_d &=&	\left[ \begin{array}{cccccccc}
									1	& 0 & 0 & 0 & 0 & 0 & 0 & 0	\\
									0	& 1 & 0 & 0 & 0 & 0 & 0 & 0	\\
									0	& 0 & 0 & 0 & 0 & 0 & 0 & 0	\\
									0	& 0 & 0 & 0 & 0 & 0 & 0 & 0	\\
									1	& 0 & 0 & 0 & 0 & 0 & 0 & 0	\\
									0	& 1 & 0 & 0 & 0 & 0 & 0 & 0	\\
									0	& 0 & 1 & 0 & 0 & 0 & 0 & 0	\\
									0	& 0 & 0 & 1 & 0 & 0 & 0 & 0	\\
									0	& 0 & 0 & 0 & 1 & 0 & 0 & 0	\\
									0	& 0 & 0 & 0 & 0 & 1 & 0 & 0	\\
									0	& 0 & 0 & 0 & 0 & 0 & 1 & 0	\\
									0	& 0 & 0 & 0 & 0 & 0 & 0 & 1	
								 \end{array}
					\right], 
\end{eqnarray*}
\begin{eqnarray*}
	X_e	&=& \left[ \begin{array}{cccccccccccc}
									0	& 0 & 0 & 0 & 0 & 0 & 1 & 0 & 0 & 1 & 0 & 0 \\
									0	& 0 & 0 & 0 & 0 & 0 & 0 & 1 & 0 & 0 & 1 & 0 \\
								 \end{array}
					\right], \\	
			& & \\
	X_f	&=& \left[ \begin{array}{cccccccccccc}
									0	& 0 & 0 & 0 & 1 & 0 & 0 & 0 & 0 & 0 & 0 & 0	\\
									0	& 0 & 0 & 0 & 0 & 1 & 0 & 0 & 0 & 0 & 0 & 0	\\
									0	& 0 & 1 & 0 & 0 & 0 & 1 & 0 & 0 & 0 & 0 & 0	\\
									0	& 0 & 0 & 1 & 0 & 0 & 0 & 1 & 0 & 0 & 0 & 0	\\
									0	& 0 & 0 & 0 & 0 & 0 & 0 & 0 & 0 & 1 & 0 & 0	\\
									0	& 0 & 0 & 0 & 0 & 0 & 0 & 0 & 0 & 0 & 1 & 0	\\
									0	& 0 & 0 & 0 & 0 & 0 & 0 & 0 & 0 & 0 & 0 & 1
								 \end{array}
					\right], \\	
			& & 	\\
	X_g	&=& \left[ \begin{array}{cccccccccccc}
									1	& 0 & 0 & 0 & 0 & 0 & 0 & 0 & 0 & 0 & 0 & 0	\\
									0	& 1 & 0 & 0 & 0 & 0 & 0 & 0 & 0 & 0 & 0 & 0	\\
									0	& 0 & 0 & 0 & 1 & 0 & 1 & 0 & 0 & 0 & 0 & 0	\\
									0	& 0 & 0 & 0 & 0 & 1 & 0 & 1 & 0 & 0 & 0 & 0	\\
									0	& 0 & 0 & 0 & 0 & 0 & 0 & 0 & 0 & 1 & 0 & 0	\\
									0	& 0 & 0 & 0 & 0 & 0 & 0 & 0 & 0 & 0 & 1 & 0	\\
									0	& 0 & 0 & 0 & 0 & 0 & 0 & 0 & 1 & 0 & 0 & 1
								 \end{array}
					\right].
\end{eqnarray*}

\noindent In particular, we see that $X_\a$ is a tree representation. 

The representation $X_\a$, as given above, is defined over every field $k$. Moreover, it is not difficult to see that $\End_{k Q} (X_\a)$ is local. Hence, the representation $X_\a$ is the unique indecomposable representation of dimension vector $\a$ over every field $k$.

Moreover, $\dim \End_{k Q} (X_\a) =9$ so that $X_\a$ is not a real Schur representation. 

\begin{theorem}
\label{main-theorem-2}
There exists no real Schur root $\b$ with the following properties:
\begin{enumerate}[(i)]
\item $X_\a \in \m^{X_\b}_{X_\b}$, \; \textrm{and}
\item $\Hom_{kQ}(X_\a,X_\b)\ne 0$\quad \textrm{or}\quad $\Hom_{kQ}(X_\b,X_\a)\ne 0$.
\end{enumerate}
\end{theorem}

\noindent If we had a sequence of real Schur roots $\b_1,\ldots,\b_n\; (n\ge 2)$ such that \linebreak $	X_\a = \sigma_{X_{\b_n}}\cdot \ldots\cdot \sigma_{X_{\b_2}} (X_{\b_1})$ then $\b_n$ would have to satisfy conditions (i) and (ii).  Note that condition (ii) merely states that $\sigma^{-1}_{X_{\b_n}}(X_\a)\ne X_\a$. Thus, once we have established the claim it is clear that $X_\a$ provides an example which answers Question $(\star )$ negatively. 

We use the rest of this section to prove the above theorem. We show that there are no real Schur roots satisfying (i). 

\begin{proof}[Proof of Theorem \ref{main-theorem-2}]

Condition (i) requires $\b < \a$ by \cite[Lemma 2]{ringel} and 
\begin{equation*}
\Ext^1_{kQ}(X_\a,X_\b) = 0 = \Ext^1_{kQ}(X_\b,X_\a),
\end{equation*}
which implies that $\langle \a,\b\rangle \ge 0$ and $\langle \b,\a \rangle \ge 0$. Hence, we start by determining the set of real roots $\b$ with the following properties:
\begin{enumerate}[(i')]
\item $\b < \a$,
\item $\langle \a,\b\rangle \ge 0$ and $\langle \b,\a \rangle \ge 0$.
\end{enumerate}

These roots are potential candidates for a reflection. Using the arguments given in \cite[Section 6]{schofield_2}, it is easy to determine the real roots $\b$ which satisfy (i') and (ii'): both conditions imply that $s_\a(\b)<0$ and, hence, if $s_\a=s_{i_1}\ldots s_{i_n}$ we get $s_\a(\b) = s_{i_1}\ldots s_{i_n}(\b)<0$ if and only if $\b=s_{i_n}\ldots s_{i_{m+1}}(e_{i_m})$ for some $m$. Thus, once we have written $s_\a$ as a product of the generators $s_i$ it is straightforward to find the real roots $\b$ satisfying (i') and (ii'). A decomposition of $s_\a$ into a product of the generators $s_i$ can be achieved as follows: if $s_i(\a)=\a'<\a$ then $s_\a=s_is_{\a'}s_i$; this gives an algorithm to find a shortest expression of $s_\a$ in terms of the $s_i$.

Applying the above algorithm to the real root $\a$, we get the following potential candidates for a reflection
\begin{eqnarray*}
	\b_1&=&(0,0,0,1,2,0,1,1), \\
	\b_2&=&(0,1,1,4,7,1,4,4),  \\
	\b_3&=&(1,0,1,4,7,1,4,4), \quad \textrm{and}	\\
	\b_4&=&(1,1,0,4,7,1,4,4).
\end{eqnarray*}

We see that $\langle \b_i,\a \rangle = 0 = \langle \a,\b_i \rangle$ for $i=2,3,4$, and hence the only reflection candidate is $\b_1$. Note that $\b_1$ is a real Schur root, and hence indeed a candidate for a reflection. However, $\b_1$ does not satisfy condition (i), that is $X_\a \notin \m^{X_{\b_1}}_{X_{\b_1}}$. Assume to the contrary that $X_\a \in \m^{X_{\b_1}}_{X_{\b_1}}$. Then $\sigma^{-1}_{X_{\b_1}} (X_\a) \in \m^{-X_{\b_1}}_{-X_{\b_1}}$, that is
\begin{eqnarray*}
\Hom_{kQ} (\sigma^{-1}_{X_{\b_1}} (X_\a),X_{\b_1}) = 0 = \Hom_{kQ} (X_{\b_1},\sigma^{-1}_{X_{\b_1}} (X_\a)).
\end{eqnarray*}
Using formula (\ref{refl-dim}) from Section \ref{universal}, we get $\g_1:=\udim \sigma^{-1}_{X_{\b_1}} (X_\a) = (1,1,1,3,2,2,2,2)$. The following diagram, however, shows that $\Hom_{kQ}(X_{\b_1},X_{\g_1})\ne 0$. The representation $X_{\g_1}$ can be constructed using the result 
\cite[Proposition A.4]{crawley} together with the same reasoning as for $X_\a$ to pass to any field $k$.

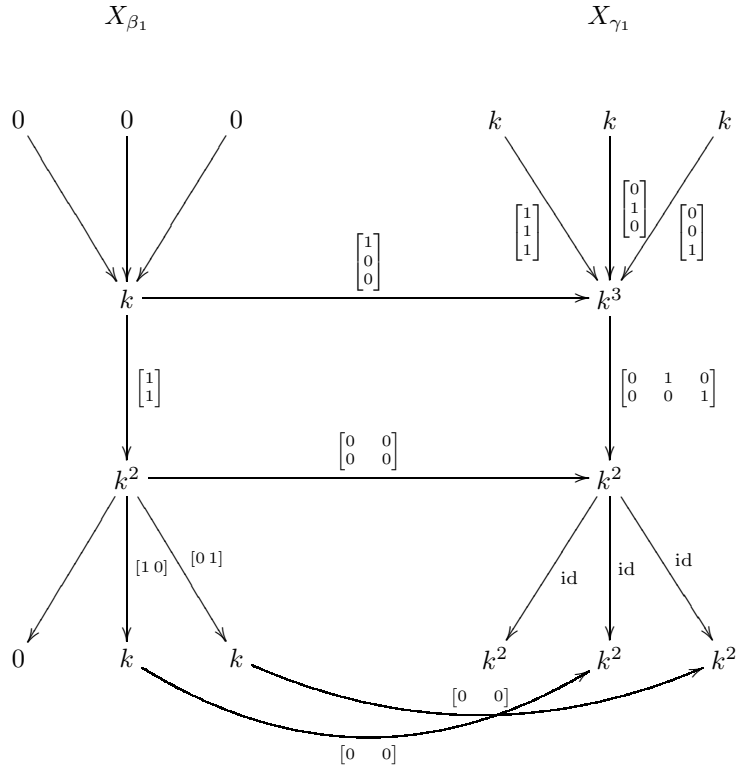
\begin{center}
$   \xymatrix{  & X_{\b_1}&&&&& X_{\g_1} &	\\
								0  \ar@<.0ex>[ddr]        &0  	  \ar@<.0ex>[dd]          		                  	     &0  \ar@<.0ex>[ddl] &&&
							  k  \ar@<.0ex>[ddr]_{\textrm{\tiny{$\begin{bmatrix}  1\\ 1\\ 1	\end{bmatrix}$}}}        &	k 
							  \ar@<.0ex>[dd]^{\textrm{\tiny{$\begin{bmatrix}  0\\ 1 \\0	\end{bmatrix}$}}}       &k \ar@<.0ex>[ddl]^{\textrm{\tiny{$\begin{bmatrix}  										0\\ 0\\ 1	\end{bmatrix}$}}}  \\
							  &&&&&&& \\
							    											 &k  	  \ar@<.0ex>[dd]^{\textrm{\tiny{$\begin{bmatrix}  1\\ 1	\end{bmatrix}$}}} 																								\ar@<.0ex>[rrrrr]^{\textrm{\tiny{$\begin{bmatrix}  1\\ 0\\ 0	\end{bmatrix}$}}}	&  &&&
				   															 &k^3  	  \ar@<.0ex>[dd]^{\textrm{\tiny{$\begin{bmatrix}  0&1&0 \\ 0&0&1	\end{bmatrix}$}}}	&	\\
								&&&&&&& \\
							 													 &k^2  	\ar@<.0ex>[ddl] \ar@<.0ex>[dd]^{\textrm{\tiny{$[1\, 0]$}}} 																															\ar@<.0ex>[ddr]^{\textrm{\tiny{$[0\, 1]$}}}  \ar@<.0ex>[rrrrr]^{\textrm{\tiny{$\begin{bmatrix}  0&0 \\0& 0	\end{bmatrix}$}}}& 	 &&&
							 													 &k^2  	\ar@<.0ex>[ddl]^{\textrm{id}}	\ar@<.0ex>[dd]^{\textrm{id}}																															\ar@<.0ex>[ddr]^{\textrm{id}} \\
								&&&&&&& \\
						    0                        &k  \ar@/^-2.5pc/[rrrrr]_{\textrm{\tiny{$\begin{bmatrix} 0&0	\end{bmatrix}$}}} &k 																							\ar@/^-1.8pc/[rrrrr]^{\textrm{\tiny{$\begin{bmatrix} 0&0	\end{bmatrix}$}}} &&& 
						    k^2                      &k^2                                      				         &k^2}$
\end{center}
This is a contradiction, and hence $X_\a\notin \m^{X_{\b_1}}_{X_{\b_1}}$ which completes the proof of the theorem and we see that, indeed, the representation $X_\a$ answers Question $(\star)$ negatively. 

\end{proof}


\begin{thebibliography}{1}
    
    \bibitem{crawley} W.W. Crawley-Boevey, `Geometry of the moment map for representations of quivers', {\em Composito	Mathematica} 126 (2001) 			257-293.
    
    
  	\bibitem{kac} V.G. Kac, `Infinite root systems, representations of graphs and invariant theory', {\em Inventiones mathematicae} 56 (1980) 					57-92. 
    
    \bibitem{ringel} C.M. Ringel, `Reflection functors for hereditary algebras', {\em J. London Math. Soc.} 21 (1980) 465-479.
    
	  \bibitem{ringel_3} C.M. Ringel, `Exceptional modules are tree modules', {\em Linear Algebra Appl.} 275/276 (1998) 471-493.
	  
    \bibitem{schofield_2} A. Schofield, `General representations of quivers', {\em Pro. London Math. Soc.} (3) 65 (1992) 46-64.
    
    \bibitem{schofield} A. Schofield, `The field of definition of a real representation of $Q$', {\em Proc. American Math. Soc.} 					       						116 (1992) 293-295.
    
	  \bibitem{wiedemann} M. Wiedemann, `Quiver representations of maximal rank type and an application to representations of a quiver with 			
	   three vertices',  {\em Bull. London Math. Soc.} 40 (2008) 479-492
	   
	  \bibitem{wiedemann-thesis} M. Wiedemann, `On real root representations of quivers', {\em PhD thesis, in preparation}
	   
\end{thebibliography}
\end{document}